\documentclass[times]{article}   	
 \pdfoutput=1
\usepackage{geometry,amsmath,amsthm} 
\usepackage{graphicx}					
\usepackage{amssymb,tikz}
\usepackage{hyperref}
\usetikzlibrary{arrows}
\usepackage{cite}

\newcommand\dyckpath[3]{

	\draw[help lines] (#1) grid +(#2,#2);
	\draw[dashed] (#1) -- +(#2,#2);
	\coordinate (prev) at (#1);
	\foreach \dir in {#3}{
		\ifnum\dir=0
		\coordinate (dep) at (1,0);
		\else
		\coordinate (dep) at (0,1);
		\fi
		\draw[line width=2pt,] (prev) -- ++(dep) coordinate (prev);
	};
}

\newcommand\partitionfr[1]{
	\coordinate (prev) at (0,0);
	\foreach \dir in {#1}{
		\draw[help lines, line width=1pt] (prev) -- +(0,1) coordinate (prev);
		\draw[help lines, line width=1pt] (prev)+(0,-1) grid +(\dir,0);
	};
}

\DeclareMathOperator{\coleg}{coleg}
\DeclareMathOperator{\area}{area}
\DeclareMathOperator{\weight}{weight}
\DeclareMathOperator{\coarm}{coarm}
\DeclareMathOperator{\Lt}{\widetilde L}
\DeclareMathOperator{\Lh}{\widehat L}
\DeclareMathOperator{\Ht}{\widetilde H}

\newtheorem{proposition}{Proposition}
\newtheorem{theorem}{Theorem}

\title{The Delta Conjecture at $q = 1$}
\author{Marino Romero\footnote{Research supported by NSF grant 9532049}\\
Department of Mathematics\\
University of California, San Diego \\
La Jolla, CA 92093\\
\href{mailto:mar007@ucsd.edu}{mar007@ucsd.edu} 
}

\begin{document}
\maketitle
\begin{abstract}

We use a weight-preserving, sign-reversing involution to find a combinatorial expansion of $\Delta_{e_k} e_n$ at $q=1$ in terms of the elementary symmetric function basis. We then use a weight-preserving bijection to prove the Delta Conjecture at $q=1$. The method of proof provides a variety of structures which can compute the inner product of $\Delta_{e_k} e_n|_{q=1}$ with any symmetric function. 
\end{abstract}
\section{Introduction}

To state our results, it will be good to start with a brief review of the various events that preceded our work. The main reason that the modified Macdonald polynomial basis was introduced by Garsia and Haiman in \cite{modifiedmac} was as a tool to prove the Macdonald $q,t$-Kostka conjectures. The idea was to prove that the Macdonald $q,t$-Kostka indexed by partitions of $n$ (with suitable modifications) give the distributions of irreducibles in a bigraded $S_n$-module. The resulting polynomials $\{\Ht_{\mu}\}_\mu$ turned out to be a remarkable symmetric function basis with many applications. In particular, they are Schur positive, which means that  we have the expansion
$$
\Ht_\mu [X; q,t] = \sum_{\lambda \vdash n} s_{\lambda}[X] \widetilde{K}_{\lambda, \mu}(q,t)
$$
with $\widetilde{K}_{\lambda, \mu}(q,t) \in \mathbb{N}[q,t]$.
This led Garsia and Haiman to  study the $S_n$-module of Diagonal Harmonics   and   conjecture a  formula for its Frobenius characteristic $DH_n[X;q,t]$ \cite{remarkable}. The expansion of this
polynomial in terms of the Modified Macdonald basis led to the introduction by Bergeron and Garsia \cite{sciencefiction} of the operator $\nabla$ 
and the birth of the expression
$$
DH_n[X;q,t] = \nabla e_n.
$$
The final proof of these identities, outlined by Claudio Procesi  in 1994, required algebraic geometrical tools and was finally carried out  by Mark Haiman \cite{n!}  around 2001.
Many years of progress in this subject led to the formulation by Haglund et. al \cite{shufflecon} of the Shuffle Conjecture which gives a combinatorial expression for $\nabla e_n$ in terms of parking functions. This conjecture was then refined to its compositional form in \cite{compos}. Soon after, Mark Haiman's work on the Hilbert Scheme attracted the attention of algebraic geometers to this subject, culminating to the formulation by Gorsky and Negut \cite{knots} of an infinite family of Shuffle Conjectures called the Rational Shuffle Conjectures, one for each co-prime pair of integers $(m,n)$, the original case corresponding to when $m= n+1$. 

Work from the 90's \cite{positivity} introduced a family of eigenoperators $\Delta_f$ of the modified Macdonald polynomial basis generalizing $\nabla$. This led to the discovery of the Schur positivity of $\Delta_{e_k} e_n$ for $1 \leq k \leq n$. This was soon followed by the formulation of the so-called Delta Conjecture \cite{deltacon} which gives a combinatorial formula for $\Delta_{e_k} e_n$, extending the classical Shuffle Conjecture. Recent papers by Erik Carlsson, Anton Mellit \cite{shuffle} then by Anton Mellit \cite{ratshuffle} give a proof of the Compositional Shuffle Conjecture and the Rational Shuffle Conjectures. However, the methods introduced in their work have yet to be seen adaptable to the Delta Conjecture which at the present moment  appears to be  of a degree of difficulty above all the previous Shuffle conjectures.

In this paper, we prove the Delta Conjecture with the parameter $q=1$. Our main tool is the use of a weight-preserving, sign-reversing involution. The symmetric function $\Delta_{e_k} e_n\big|_{q=1}$ is expanded into a formal power series in $t$ involving an infinite number of terms. The sign-reversing involution is used to cancel all the negative signs and output a finite number of fixed points yielding a polynomial in $t$ with positive integer coefficients. The main difficulty was  the construction of a suitable family of combinatorial objects to represent this formal series.   Our choice of objects caused 
the weight-preserving property of the involution to encode 
the needed combinatorial properties of the fixed points. 
More precisely,  properties whose defining characteristics are closely related to the combinatorial side of the conjecture. This leads to a bijection between the predicted combinatorial side and the symmetric function side of the Delta Conjecture at $q=1$. 

By these methods, we are able to compute other cases of the Delta operator at $q=1$, such as $\Delta_{h_k} e_n$. However, the importance of our work lies in introducing a new way, namely the sign-reversing involution method, for proving the Delta Conjecture and, in particular, the Shuffle Conjecture. Using similar steps, the rational function $\Delta_{e_k} e_n$ may be expressed as a formal power series in $q$ and $t$ given by a certain set of objects. In principle, one should be able to find a sign-reversing involution of these objects, providing a combinatorial expansion of $\Delta_{e_k} e_n$.

Our main result is a proof of the identity 
$$
\Delta_{e_k} e_n \Big|_{q=1}   = \sum_{\lambda \vdash n} e_\lambda \sum_{s \in M^{\lambda}_k} t^{\rho(s)}, \eqno(1)
$$
where $M^{\lambda}_k$ is the set of all sequences 
$$
s=((a_1,b_1), (a_2,b_2), \dots, (a_{k+1}, b_{k+1}))
$$ 
composed of nonnegative integers satisfying
\begin{enumerate}
\item $a_1 = 0$,
\item $(b_1,\dots, b_{k+1})$ is a rearrangement of $\lambda_1,\dots, \lambda_{\ell(\lambda)}, 0^{k+1-\ell(\lambda)}$ (for $\ell(\lambda) \leq k+1$),
\item $a_{i+1} < a_i + b_i$ and
\item  $\rho(s) = a_1+ \cdots+ a_{k+1}$.
\end{enumerate}

Let $D_n$ be the collection of Dyck paths in the $n \times n$ square. That is, $D \in D_n$ is a lattice path from $(0,0)$ to $(n,n)$  consisting of North and East steps which stay weakly above the main diagonal $y=x$. Reading from bottom to top, let $\alpha_i$ be the number of lattice cells in the $i^{\text{th}}$ row of $D$ found between the path and the main diagonal. We call $(\alpha_1,\dots, \alpha_n)$ the area sequence of $D$, and let $\area(D)$ be the sum $\alpha_1+ \cdots + \alpha_n$. We denote by $\lambda(D)$ the partition of $n$ which is made by the lengths of vertical segments in $D$. 
To exhibit these definitions, we provide the picture for the Dyck path whose area sequence is  $(0,1,2,2,1,2,1,2)$:

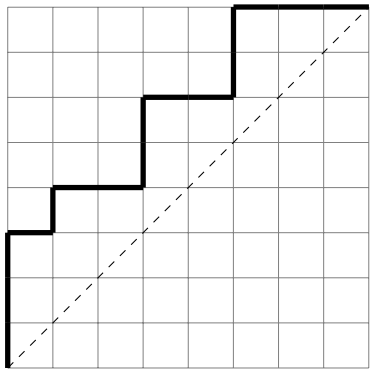
\begin{align*}
\begin{tikzpicture}[scale =.6] 
\dyckpath{0,0}{8}{1,1,1,0,1,0,0,1,1,0,0,1,1,0,0,0};
\end{tikzpicture}.
\end{align*}
In this case, $\lambda(D) = (3,2,2,1)$ since the vertical segments have lengths $3,1,2,2$.

The Delta Conjecture \cite{deltacon}, at $q=1$, may be restated as
$$
\Delta_{e_k} e_n\Big|_{q=1} = \sum_{D\in D_n} t^{\area(D)} H_{n-k}(D) e_{\lambda(D)},
$$
where if $D$ has area sequence $(\alpha_1,\dots, \alpha_n)$, then 
$$
H_{n-k}(D) = (1+w)\prod_{\alpha_{i-1}  = \alpha_i -1 } \left(1+\frac{w}{t^{\alpha_i}}\right)\Big|_{w^{n-k}}.
$$

The identity in (1) proves the Delta Conjecture at $q=1$ since we will bijectively show
$$
\sum_{s \in M_k^{\lambda}} t^{\rho(s)} = \sum_{\substack{D \in D_n \\ \lambda(D) = \lambda}} t^{\area(D)} H_{n-k}(D). \eqno(2)
$$

The order of our exposition is as follows: We begin with some definitions and manipulate a symmetric function expansion of $\Delta_{e_k}e_n$. It will become essential to use an addition formula for the forgotten basis of symmetric functions $\{f_{\mu}\}_{\mu}$ and include for completeness a section with these computations. The expansion of $\Delta_{e_k}e_n$ appears as a formal power series with alternating signs, but we find (1) with a sign-reversing involution. We will finish with a bijection, establishing (2).  

\section{Definitions and Formulations}

For readers who are not familiar with plethystic notation and other symmetric function ammenities, we refer the reader to \cite{plethysm}.
Following Macdonald's section on Integral Forms on page 364 of \cite{macdonald}, we have the equality
$$
J_{\mu}[X;1,t] = \prod_{i=1}^{\ell(\mu')} (t:t)_{\mu'_i} e_{\mu'_i}[X],
$$
where
$$
(q:t)_k = (1-q) (1- q t) \cdots (1- q t^{k-1}).
$$

The modified Macdonald polynomials of \cite{modifiedmac} arise from $J_{\mu}$ by setting
$$
\tilde{H}_{\mu} [X;q,t] = t^{n(\mu)}J_{\mu}\left[\frac{X}{1-1/t};q, 1/t\right],
$$ 

where $n(\mu) = \sum_{i=1}^{\ell(\mu)} (i-1) \mu_i$. Thus, we have

\begin{align*}
\tilde{H}_{\mu} [X;1,t] &=  t^{n(\mu)}J_{\mu}\left[\frac{X}{1-1/t};1 , 1/t\right] 
 = t^{n(\mu)} \prod_{i=1}^{\ell(\mu')} (1/t:1/t)_{\mu'_i} e_{\mu'_i}\left[\frac{X}{1-1/t}\right]\\
&=  t^{n(\mu)} \prod_{i=1}^{\ell(\mu')} t^{-(1+2+\cdots +\mu'_i)} (-1)^{\mu'_i} (t:t)_{\mu'_i} e_{\mu'_i}\left[-t \frac{X}{1-t}\right]\\
&=  t^{n(\mu)} \prod_{i=1}^{\ell(\mu')} t^{-(1+2+\cdots +\mu'_i)} (-1)^{\mu'_i} (t:t)_{\mu'_i} (-t)^{\mu'_i}h_{\mu'_i}\left[\frac{X}{1-t}\right]\\
&= t^{n(\mu)} \prod_{i=1}^{\ell(\mu')} t^{-\binom{\mu'_i}{2}} (t:t)_{\mu'_i} h_{\mu'_i}\left[\frac{X}{1-t}\right]
\\\end{align*} giving the equality \begin{align*}
\tilde{H}_{\mu} [X;1,t] &= \prod_{i=1}^{\ell(\mu')} (t:t)_{\mu'_i} h_{\mu'_i}\left[\frac{X}{1-t}\right].
\end{align*}
\\

The Delta operator, first seen in \cite{positivity} and discussed by Haglund \cite{schroder} and Wilson, 
Haglund, Remmel  \cite{deltacon},  is defined for any symmetric function $f$ as the eigen-operator on the modified Macdonald polynomials given by
$$
\Delta_f \tilde{H}_\mu[X;q,t] = f[B_\mu(q,t)]\tilde{H}_\mu[X;q,t],
$$
where $B_\mu(q,t) = \sum_{c \in \mu} t^{\coleg_\mu(c)} q^{\coarm_\mu(c)}$. In particular, it follows that 
$$
\Delta_f h_{\mu'}\left[\frac{X}{1-t} \right] \Big|_{q=1}  =  f[B_\mu(1,t)]h_{\mu'}\left[\frac{X}{1-t} \right]\\
$$
Switching the roles of $\mu$ and $\mu'$, and using that 
$$B_{\mu'}(1,t) = \sum_{c \in \mu'} t^{\coleg_{\mu'}(c)} = \sum_{i=1}^{\ell(\mu)} [\mu_i]_t,$$
we get
$$
\Delta_f h_{\mu}\left[\frac{X}{1-t}\right] \Big|_{q=1} = f\left[ \sum_{i=1}^{\ell(\mu)} [\mu_i]_t\right] h_{\mu}\left[\frac{X}{1-t}\right].
$$

This means that $\Delta_f e_n$ can be computed from an expansion of $e_n$ in terms of $\left\{ h_{\mu}\left[\frac{X}{1-t}\right]\right\}_{\mu}.$
To this end, we note that for any expressions $X,Y$,
$$
e_n[XY] = h_n[X (-\epsilon Y)] \Big|_{\epsilon = -1} = \sum_{\mu \vdash n} h_\mu[X] m_\mu[-\epsilon Y] \Big|_{\epsilon = -1}= \sum_{\mu \vdash n} h_\mu[X] f_\mu[ Y] ,
$$
where $f_{\mu}$ is the forgotten symmetric function indexed by $\mu$. Thus
\begin{align*}
e_n[X] = e_n\left[\frac{X}{1-t} (1-t)\right] & = \sum_{\mu \vdash n} h_\mu\left[ \frac{X}{1-t}\right] f_\mu[ 1-t]. 
\end{align*}

For any symmetric function $F \in \Lambda^n$, Haglund's Identity \cite{schroder} gives the remarkable equality
$$
\langle \Delta_{e_k}e_n , F \rangle = 
\langle \Delta_{\omega(F)} e_{k+1}, s_{k+1} \rangle.
$$

It therefore follows that the coefficient of $e_\lambda$  in the expansion of $\Delta_{e_k}e_n $ at $q=1$ equals
$$
\langle \Delta_{e_k}e_n , f_{\lambda} \rangle \Big|_{q=1} = 
\langle \Delta_{\omega(f_{\lambda})} e_{k+1}, s_{k+1} \rangle\Big|_{q=1} 
=
\langle \Delta_{m_{\lambda}} e_{k+1}, s_{k+1} \rangle\Big|_{q=1} 
$$
$$
= \langle 
 \sum_{\mu \vdash k+1} h_\mu\left[ \frac{X}{1-t}\right] f_\mu[ 1-t] m_{\lambda}\left[\sum [\mu_i]_t \right]
 , s_{k+1}
\rangle
$$

However, with an application of the Cauchy formula together with $\langle h_{\lambda}, h_{k+1} \rangle = 1$, we get
$$
\langle h_{\mu}\left[ \frac{X}{1-t} \right], s_{k+1} \rangle =  h_{\mu}\left[ \frac{1}{1-t} \right].
$$
This gives the beginning stage of our desired coefficient:
$$
 \langle \Delta_{e_k}e_n , f_{\lambda} \rangle \Big|_{q=1} = \sum_{\mu \vdash k+1} h_{\mu} \left[ \frac{1}{1-t} \right] f_\mu [1-t] m_{\lambda} \left[\sum [\mu_i]_t \right]
$$
\\

We note here that if $\ell(\lambda) > k+1$, then $m_{\lambda} \left[\sum [\mu_i]_t \right] = 0$. So we safely assume $\lambda$ is a partition whose length is no more than $k+1$.  

\section{An addition formula for the Forgotten}

We wish to find the value of $f_{\mu}[1-t]$ for $\mu \vdash k+1$. For the remainder, we avoid repeatedly writing $k+1$ and we set $m = k+1$. One can follow the results from Remmel and E{\~ g}ecio{\~ g}lu \cite{bricktabloids}, but for completion, we give a separate account.

Suppose we have expressions $X,Y,$ and $Z$. Then 
$$
e_{m}[X(Y+Z)] =  \sum_{\mu \vdash m} h_\mu[X] f_\mu[ Y + Z].
$$
On the other hand,
\begin{align*}
e_m[XY+XZ]  = && \sum_{i=0}^m e_i[XY] e_{m-i}[YZ] 
 = && \sum_{i=0}^m \sum_{\alpha \vdash  i}\sum_{\beta \vdash m-i} h_\alpha[X] f_\alpha[Y]  h_\beta [X] f_\beta[Z].
\end{align*}
For $ h_\alpha[X] h_\beta[X] = h_{\mu}[X]$, we require that the parts of $\alpha$ and $\beta$ together create $\mu$. Therefore, allowing the empty partition and letting $(\alpha , \beta)$ be the partition obtained from the union of parts from $\alpha$ and $\beta$, we can write
$$
f_{\mu}[Y+Z] = \sum_{(\alpha,\beta) = \mu} f_{\alpha}[Y] f_{\beta}[Z].
$$

We note that 
$$
e_m[X (1) ] = \sum_{\mu \vdash m} h_{\mu}[X] f_{\mu}[1],
$$
meaning that $f_{\mu}[1]$ is the coefficient of $h_{\mu}$ in the Jacobi-Trudi determinant:
$$
\det \begin{bmatrix}
h_1 && h_2 && \cdots && h_m\\
1 && h_1 && \cdots && h_{m-1}\\
\vdots && \ddots&& \ddots && \vdots\\
0 && \cdots && 1&& h_1
\end{bmatrix}.
$$
Though this may be computed using Brick Tabloids \cite{bricktabloids}, we give chance to another amazing way of computing this coefficient:

In proving a formula for the number of partitions in a given shape Sergel Leven and Amdeberhan \cite{tewodros} notice that for any matrix whose elements below the sub-diagonal are $0$, every term in the determinant is uniquely determined by a selection of elements along the sub-diagonal. For example, suppose we select the $1's$ in rows $2,3,5,6,7$ in expanding $e_{8}$. We get the following picture:

\begin{align*}
\begin{tikzpicture}[scale =.5] 
	\draw [line width=0.25mm ] (-5.7,.8) -- (-1.7,.8)-- (-1.7,3.5) -- (-5.7 , 3.5)--(-5.7,.8);
	\draw [line width=0.25mm ] (-1.3,-2.5) -- (4.2,-2.5)-- (4.2,.9) -- (-1.3 , .9)--(-1.3,-2.5);
	\draw [line width=0.25mm ] (4.4,-3.3) -- (5.5,-3.3)-- (5.5,-2.5) -- (4.4 ,-2.5)--(4.4,-3.3);
\node at (0,0) {{$
\begin{bmatrix}
&h_1 & h_2 & h_3 & h_4 & h_5 & h_6 & h_7 & h_8 & \\ 
&\underline{1} & h_1 & h_2  & h_3 & h_4 & h_5 & h_6 & h_7&\\
&0 & \underline{1} & h_1  & h_2 & h_3 & h_4 & h_5 & h_6&\\
&0 & 0 & 1 & h_1  & h_2 & h_3 & h_4 & h_5&\\
&0 & 0 & 0 & \underline{1} & h_1  & h_2 & h_3 & h_4&\\
&0 & 0 & 0 & 0 & \underline{1} & h_1  & h_2 & h_3&\\
&0 & 0 & 0 & 0 & 0 & \underline{1}  & h_1 & h_2&\\
&0 & 0 & 0 & 0 & 0 & 0 &1 &  h_1&
\end{bmatrix}
$}};
\end{tikzpicture} 
\end{align*}

Proceeding from left to right, we first see that columns $1$ and $2$ have an entry selected. Column $3$ has no $1$ selected, so it is then forced to take $h_3$. Now all first three rows have an element selected. The next column with no $1$ selected is $7$. The three previous $1$'s are selected, so we must take $h_4$. Lastly, we take $h_1$ at the bottom. Each block forms a cycle in the permutation that selects these elements in the determinant.
We can therefore break up the determinant into squares of sizes $\mu_1,\dots, \mu_{\ell(\mu)}$ along the main diagonal, each contributing a sign of $(-1)^{\mu_i - 1}$ and a factor of $h_{\mu_i} $.

From this, it follows that we can compute that $f_{\mu}[1]$ equals 
$$
f_{\mu}[1] = (-1)^{m- \ell(\mu)} |R(\mu)|
$$
where $R(\mu)$ is the collection of rearrangements of the parts $\mu_1,\dots, \mu_{\ell(\mu)}$. We apply this to the summation expansion

\begin{align*}
f_{\mu}[1-t] & = \sum_{(\alpha,\beta) = \mu} f_{\alpha}[1] f_{\beta}[-t] \\ 
& = \sum_{(\alpha,\beta) = \mu} (-t)^{\beta} f_{\alpha}[1] m_{\beta}[1].
\end{align*}
Since $m_{\beta}[1]$ is $1$ if $\beta$ has only one part and $0$ otherwise, we can write
$$
f_{\mu}[1-t] =f_{\mu}[1] + \sum_{i \in \{\mu_1,\dots, \mu_n\}} (-t)^i f_{\mu - (i)}[1],
$$
where $\mu - (i)$ is the partition obtained from $\mu$ by removing a part of size $i$. We will denote the inclusion in the sum as $i \in \{\mu\}$. Our observation from the Jacobi-Trudi determinant then gives
$$
f_{\mu}[1-t] =(-1)^{m- \ell(\mu)}|R(\mu)|+ \sum_{i \in \{\mu\}} t^i (-1)^i (-1)^{(m-i) - (\ell(\mu)-1)}|R(\mu-(i))|$$
$$
=(-1)^{m- \ell(\mu)} \left( |R(\mu)| - \sum_{i \in \{\mu\}} t^i |R(\mu-(i))| \right).
$$

\section{$h_{\mu}\left[\frac{1}{1-t}\right] f_\mu[1-t]$}

Since 
$$
h_{r}\left[\frac{1}{1-t}\right]  = \frac{1}{(1-t)(1-t^2) \cdots (1-t^r)},
$$ 
we can think of this factor, which we denote by $G_r$, to be the generating function of all partitions whose largest row is less than or equal to $r$. In other words,
$$
G_r = h_r\left[ \frac{1}{1-t} \right] = \sum_{\lambda: \lambda_1 \leq r} t^{|\lambda|}.
$$
Using the expression from the previous section, we get
$$
h_\mu\left[\frac{1}{1-t} \right] f_{\mu}[1-t] = (-1)^{m-\ell(\mu)}G_{\mu_1}\cdots G_{\mu_{\ell(\mu)}} \left( |R(\mu)|
- \sum_{i \in \{\mu\}}t^{i} |R(\mu- (i))| \right).
$$
\begin{proposition}
$$h_\mu\left[\frac{1}{1-t} \right] f_{\mu}[1-t]  = (-1)^{m-\ell(\mu)} \sum_{i \in \{\mu\}}
G_{i-1} \frac{ G_{\mu_1} \cdots G_{\mu_{\ell(\mu)}}}{G_i} |R(\mu -(i))|$$
\begin{proof}

To begin introducing the combinatorial objects which we will use to expand the main symmetric function identity, we will now give an over-embellished method of describing this particular sum.
Let $L_\mu$ denote the objects formed by 
\begin{enumerate}
\item selecting a rearrangement $(\mu_{\alpha_1}, \dots, \mu_{\alpha_{{\ell(\mu)}}}) $ of $\mu_1,\dots, \mu_{\ell(\mu)}$, 
\item choosing a sequence $(\nu^{1},\dots, \nu^{\ell(\mu)})$ of partitions satisfying $\nu^i_1 \leq \mu_{\alpha_i}$, and
\item drawing the sequence of Ferrer's diagrams $(\nu^{1},\dots, \nu^{\ell(\mu)})$, where $\nu^i$ is drawn above a row of size $\mu_{\alpha_i}$.
\end{enumerate}

 For example, given the partition $\mu= (3,2,2,1)$ and rearrangement $(2,1,3,2)$ we could find in $L_{\mu}$ the sequence $((2,1,1,1), (0), (2,1,1),(1,1,1,1))$ which we depict by
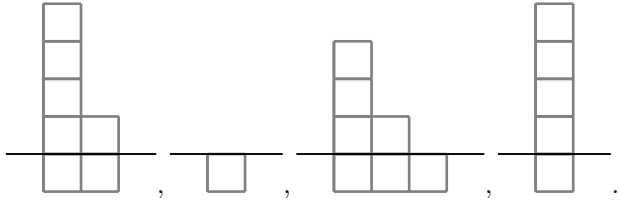
\begin{align*}
\begin{tikzpicture}[scale =.5] 
	\partitionfr{2,2,1,1,1}
	\draw [line width=0.25mm ] (-1,1) -- (3,1);
\end{tikzpicture},
\begin{tikzpicture}[scale =.5] 
	\partitionfr{1}
	\draw [line width=0.25mm ] (-1,1) -- (2,1);
\end{tikzpicture},
\begin{tikzpicture}[scale =.5] 
	\partitionfr{3,2,1,1}
	\draw [line width=0.25mm ] (-1,1) -- (4,1);
\end{tikzpicture},
\begin{tikzpicture}[scale =.5] 
	\partitionfr{1,1,1,1,1}
	\draw [line width=0.25mm ] (-1,1) -- (2,1);
\end{tikzpicture}.
\end{align*}

For $T \in L_{\mu} $ with partitions $\nu^{1},\dots, \nu^{{\ell(\mu)}}$, let $w(T) = |\nu^{1}|+\cdots+ |\nu^{{\ell(\mu)}}|$.  For instance, the weight of our example is $5+0+4+4 = 13$. Then
$$
G_{\mu_1} \cdots G_{\mu_{\ell(\mu)}} |R(\mu)| = \sum_{T \in L_{\mu}} t^{w(T)} = \sum_{i \in \{\mu\}} G_i \sum_{S \in L_{\mu-(i)}} t^{w(S)},
$$
where the last equality separates the sum in terms of which part of $\mu$ begins the rearrangement ($\mu_{\alpha_1} = i$).
On the other hand 
$$
\sum_{i\in \{\mu\}} t^i G_{\mu_1} \cdots G_{\mu_{\ell(\mu)}}  |R(\mu-(i))|
= \sum_{i \in \{\mu\}} t^i G_i  \sum_{S \in L_{\mu-(i)}} t^{w(S)}.
$$
But $t^i G_i$ can be thought of as the generating function of all partitions whose biggest part equals $i$. So we have $G_i - t^i G_i = G_{i-1}$, giving the proposition:
$$h_\mu\left[\frac{1}{1-t} \right] f_{\mu}[1-t]  = (-1)^{m-\ell(\mu)} \sum_{i \in \{\mu\}}
G_{i-1}\sum_{S \in L_{\mu - (i)}} t^{w(S)}
$$
\end{proof}
\end{proposition}

Let $\Lh_{\mu}$ be the subset of objects from $L_{\mu}$ with the restriction that the first partition in our sequence, say $\nu^{1}$ corresponding to the part $\mu_{\alpha_1}$, satisfies $\nu^{1}_1 \leq \mu_{\alpha_1}-1$ instead of $\mu_{\alpha_1}$. Then it follows from the proposition that
$$
h_\mu\left[\frac{1}{1-t} \right] f_{\mu}[1-t] = (-1)^{m- \ell(\mu)} \sum_{T \in \Lh_{\mu}} t^{w(T)}.
$$

\section{The summation terms}
We now describe the product
$$h_{\mu} \left[ \frac{1}{1-t} \right] f_\mu [1-t] m_{\lambda} \left[\sum [\mu_i]_t \right]$$ 
by labeling the objects in $\Lh_{\mu}$. Let $\Lt_{\mu}$ be the set of objects formed in the following way:
\begin{enumerate}
\item Select an element in $\Lh_{\mu}$ (recall that $\nu^1_1 \leq \mu_{\alpha_1}-1$).
\item Place each of $\lambda_1,\dots, \lambda_{\ell(\lambda)}$ in a cell of $\mu_1,\dots,\mu_{\ell(\mu)}$ such that each cell contains at most $1$ entry. Fill the remaining cells with a $0$.  
\item If $\lambda_i$ is placed in a cell with $j$ cells to the left, then we say that it contributes $j \cdot \lambda_i$.
\end{enumerate}
For $T\in\Lt_\mu$, we define the weight $p(T)$ by
$$
p(T) = w(T)+\sum_{\lambda_i} (\text{the contribution of $\lambda_i$}).
$$

For example, let us select the partition $\mu = (4,2,2,1)$ and $\lambda=(5,3,3,2,1,1)$.
Then we can select a rearrangement of $\mu$: $(2,1,4,2)$. Then select partitions 
$ (\nu^1, \nu^2, \nu^3, \nu^4)$ so that the largest part of $\nu^1$ is at most $2-1$, the largest part of $\nu^2$ is at most $1$,the largest part in $\nu^3$ is at most $4$ and the largest part in $\nu^4$ is at most $2$. Here we choose 
$$
((1,1,1),(1),(3,2,2,1),(2,2,1)).
$$ 
We then place the parts of $\lambda$ in the cells corresponding to $\mu_1,\dots,\mu_{\ell(\mu)}$. One such selection is given by
\begin{align*}
\begin{tikzpicture}[scale =.5] 
	\partitionfr{2,1,1,1}
	\draw [line width=0.25mm ] (-1,1) -- (3,1);
	\node at (.5,.5) {{$0$}};
	\node at (1.5,.5) {{$2$}};
	\node at (-2.5,2) {{$S =$}};
\end{tikzpicture},
\begin{tikzpicture}[scale =.5] 
	\partitionfr{1,1}
	\draw [line width=0.25mm ] (-1,1) -- (2,1);
	\node at (.5,.5) {{$5$}};
\end{tikzpicture},
\begin{tikzpicture}[scale =.5] 
	\partitionfr{4,3,2,2,1}
	\draw [line width=0.25mm ] (-1,1) -- (5,1);
	\node at (.5,.5) {{$1$}};
	\node at (1.5,.5) {{$3$}};
	\node at (2.5,.5) {{$0$}};
	\node at (3.5,.5) {{$3$}};
\end{tikzpicture},
\begin{tikzpicture}[scale =.5] 
	\partitionfr{2,2,2,1}
	\draw [line width=0.25mm ] (-1,1) -- (3,1);
	\node at (.5,.5) {{$0$}};
	\node at (1.5,.5) {{$1$}};
\end{tikzpicture}.
\end{align*}

We get $w(S) = 17$ since the sum of sizes of the partitions is $17$. The contributions from the fillings (from left to right respectively) is given by $$0 \cdot (0) +1\cdot(2),~ 0 \cdot(5), ~0\cdot(1)+1\cdot(3)+2\cdot(0)+3\cdot(3), ~0\cdot(0)+1\cdot(1)$$
or $
2,0,12,1
$. Therefore $p(S) = 17+(2+0+12+1) = 32$. We say that the sign associated to $T \in \Lt_{\mu}$ is $(-1)^{m-\ell(\mu)}$. So in this case, the sign is given by $(-1)^{9-4 } = -1$.

We then get, for $\mu \vdash m$,
\begin{proposition}
$$h_{\mu} \left[ \frac{1}{1-t} \right] f_\mu [1-t] m_{\lambda} \left[\sum [\mu_i]_t \right]
= (-1)^{m-\ell(\mu)} \sum_{T \in \Lt_{\mu}} t^{p(T)}.$$
\begin{proof}
We have
$$
m_{\lambda}(x_1,\dots,x_m) = \sum_{(\lambda_{i_1},\dots, \lambda_{i_m})} x_1^{\lambda_{i_1}} \cdots x_m^{\lambda_{i_m}},
$$
where the sum is over all rearrangements $(\lambda_{i_1},\dots, \lambda_{i_m})$ of $\lambda_1,\dots, \lambda_m = \lambda,0^{m-\ell(\lambda)}$.
We can interpret this as labeling $x_1,\dots, x_m$ with $\lambda_1,\dots, \lambda_m$. Now we substitute $(x_1,x_2, \cdots , x_{\mu_1}) = (1,t,\dots, t^{{\mu_1}-1})$, then $(x_{\mu_1+1},\dots, x_{\mu_1+\mu_2}) = (1,t,\dots, t^{\mu_2-1})$, and so on. For example, suppose we are computing $m_{3,2,1,1}[[3]_t+[2]_t+[2]_t]$. Then a monomial is given by a rearrangement of $3,2,1,1,0,0,0$. Selecting the rearrangement $(0,1,0,2,3,0,1)$ corresponds to taking the monomial 
$$x_1^0x_2^1x_3^0x_4^2x_5^3x_6^0 x_7^1.$$ Writing this rearrangement in the cells of $\mu_1=3, \mu_2=2,$ and $\mu_3=2 $, we get the filling
\begin{align*}
\begin{tikzpicture}[scale =.5] 
	\partitionfr{3}
	\node at (.5,.5) {{$0$}};
	\node at (1.5,.5) {{$1$}};
	\node at (2.5,.5) {{$0$}};
	\node at (.5,-.5) {\scriptsize{$1$}};
	\node at (1.5,-.5) {\scriptsize{$t$}};
	\node at (2.5,-.5) {\scriptsize{$t^2$}};
\end{tikzpicture}&&
\begin{tikzpicture}[scale =.5] 
	\partitionfr{2}
	\node at (.5,.5) {{$2$}};
	\node at (1.5,.5) {{$3$}};
	\node at (.5,-.5) {\scriptsize{$1$}};
	\node at (1.5,-.5) {\scriptsize{$t$}};
\end{tikzpicture}&&
\begin{tikzpicture}[scale =.5] 
	\partitionfr{2}
	\node at (.5,.5) {{$0$}};
	\node at (1.5,.5) {{$1$}};
	\node at (.5,-.5) {\scriptsize{$1$}};
	\node at (1.5,-.5) {\scriptsize{$t$}};
\end{tikzpicture}
\end{align*}

We wrote under each cell the term which it represents in $[\mu_1]_t +[\mu_2]_t +[\mu_3]_t $. Reading from left to right, we get the term
$$
1^0 t^1 (t^2)^0 1^2 t^3 1^0 t^1.
$$

This means 
$
m_{\lambda}[\sum[\mu_i]_t]
$
can be interpreted as the sum over fillings of $\mu_1, \dots, \mu_{\ell(\mu)}$ by $\lambda_1,\dots, \lambda_{m}$.  If $\mu_i$ is filled by $(\lambda_{i_1},\dots, \lambda_{i_{\mu_1}})$, the filling of $\mu_i$ will contribute $$(1^{\lambda_{i_1}} (t)^{\lambda_{i_2}} \cdots (t^{\mu_i-1})^{\lambda_{i_{\mu_i}}})$$ to the term in the sum. This corresponds to saying that if $\lambda_{i}$ is placed in a cell with $j$ cells to the left, then it contributes to the term by $(t^j)^{\lambda_i} $.
\\

The sum over the contributions of $\lambda_i$ gives exactly the power of $t$ contributed by a monomial in the expansion of $m_{\lambda}$.
\end{proof}
\end{proposition}
We have the essential formulation
$$
\langle \Delta_{e_k} e_n, f_{\lambda} \rangle \Big|_{q=1}
= \sum_{\mu \vdash k+1} (-1)^{k+1- \ell(\mu)} \sum_{T \in \Lt_{\mu}} t^{p(T)}.
$$

We now use a sign-reversing involution to reduce this series to a positive polynomial. 

\section{The involution}
The involution follows from the following observation:

Given a partition $\nu$ over a row of size $\mu_i > 1$ and a labeling of $\mu_i$ whose last cell is labeled by $c$, we can create a pair of partitions $\nu^1, \nu^2$ over rows of sizes $1$ and $\mu_i -1$ respectively, by letting the last column in $\nu$ be $\nu^1$, and letting $\nu^2$ be the partition obtained from $\nu$ by removing the last column and adding $c$ rows of length $\mu_i -1$. We then add $c$ to the cell under $\nu^1$. To better describe this long-winded process, we will walk through a particular example. Suppose we have the partition $(4,4,2,1)$ above a row of length $\mu_i = 4$. And suppose we have filled $\mu_i$ with $(0,1,0,2)$ so that the picture looks as follows:\\

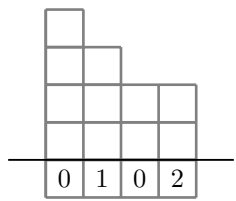
\begin{align*}
\begin{tikzpicture}[scale =.5] 
	\partitionfr{4,4,4,2,1}
	\draw [line width=0.25mm ] (-1,1) -- (5,1);
	\node at (.5,.5) {{$0$}};
	\node at (1.5,.5) {{$1$}};
	\node at (2.5,.5) {{$0$}};
	\node at (3.5,.5) {{$2$}};
\end{tikzpicture}
\end{align*}

In this case $c = 2$. We then make the last column ($\nu^1$ over $c$) the first part of our image and add $c = 2$ rows of length $3$ to the partition over $(0,1,0)$:
\\
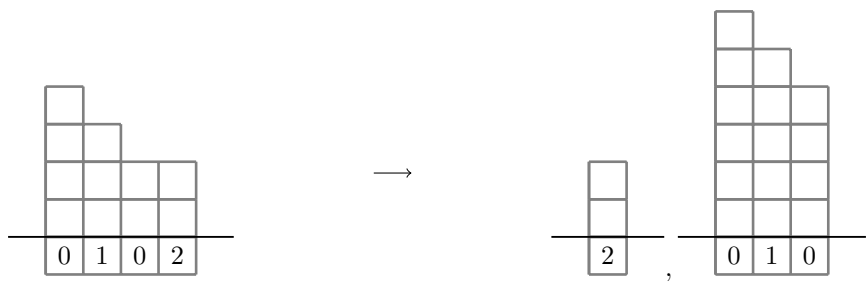
\begin{align*}
\begin{tikzpicture}[scale =.5] 
	\partitionfr{4,4,4,2,1}
	\draw [line width=0.25mm ] (-1,1) -- (5,1);
	\node at (.5,.5) {{$0$}};
	\node at (1.5,.5) {{$1$}};
	\node at (2.5,.5) {{$0$}};
	\node at (3.5,.5) {{$2$}};
\end{tikzpicture}
&&
\begin{tikzpicture}[scale =.5] 
\node at (.5,.5) {{}};
\draw [-> ] (0,3) -- (1,3);
\end{tikzpicture}
&&
\begin{tikzpicture}[scale =.5] 
	\partitionfr{1,1,1}
	\draw [line width=0.25mm ] (-1,1) -- (2,1);
	\node at (.5,.5) {{$2$}};
\end{tikzpicture},
\begin{tikzpicture}[scale =.5] 
	\partitionfr{3,3,3,3,3,2,1}
	\draw [line width=0.25mm ] (-1,1) -- (4,1);
	\node at (.5,.5) {{$0$}};
	\node at (1.5,.5) {{$1$}};
	\node at (2.5,.5) {{$0$}};
\end{tikzpicture}
\end{align*}
\\
This operation preserves the weight of the object since the contribution of $c$ in the pre-image is equal to $c \cdot( \mu_i-1)$, which is given by the $c$ rows of length $\mu_i-1$ that were added in the second partition.
We note that given a column $\nu^1$ over a cell filled with $c$ and partition $\nu^2$ over a row of length $\mu_{i}-1$,
we can combine the two, inverting this procedure, provided that we can remove $c$ full rows from $\nu_2$ and add a column of size $|\nu^1|$ on the right. This condition is described by the inequality
$$
|\nu^1| \leq (\text{ the number of rows in $\nu^2$ of size $\mu_{i}-1$}) -c.
$$
We instead write the inequality as
$$
c + |\nu^1| \leq (\text{ the number of rows in $\nu^2$ of size $\mu_{i}-1$})
$$
For example,
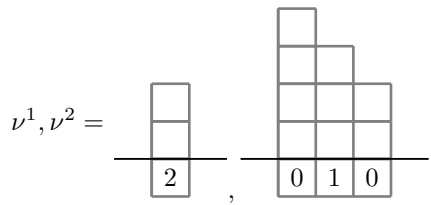
\begin{align*}
\begin{tikzpicture}[scale =.5] 
	\partitionfr{1,1,1}
	\draw [line width=0.25mm ] (-1,1) -- (2,1);
	\node at (.5,.5) {{$2$}};
	\node at (-2.5,2) {{$\nu^1,\nu^2 =$}};
\end{tikzpicture},
\begin{tikzpicture}[scale =.5] 
	\partitionfr{3,3,3,2,1}
	\draw [line width=0.25mm ] (-1,1) -- (4,1);
	\node at (.5,.5) {{$0$}};
	\node at (1.5,.5) {{$1$}};
	\node at (2.5,.5) {{$0$}};
\end{tikzpicture}
\end{align*}
cannot be combined since the number of rows in $\nu^2$ of length $3$ is $2$, and $c+|\nu^1| = 2 + 2 = 4 >2 $. However, we can expand the second partition because the length of the row underneath it is greater than $1$. We get
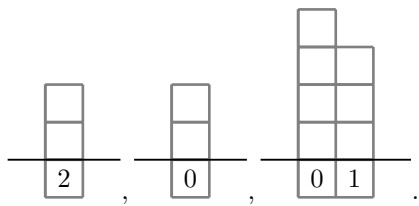
\begin{align*}
\begin{tikzpicture}[scale =.5] 
	\partitionfr{1,1,1}
	\draw [line width=0.25mm ] (-1,1) -- (2,1);
	\node at (.5,.5) {{$2$}};
\end{tikzpicture},
\begin{tikzpicture}[scale =.5] 
	\partitionfr{1,1,1}
	\draw [line width=0.25mm ] (-1,1) -- (2,1);
	\node at (.5,.5) {{$0$}};
\end{tikzpicture},
\begin{tikzpicture}[scale =.5] 
	\partitionfr{2,2,2,2,1}
	\draw [line width=0.25mm ] (-1,1) -- (3,1);
	\node at (.5,.5) {{$0$}};
	\node at (1.5,.5) {{$1$}};
\end{tikzpicture}.
\end{align*}
The crucial observation is that the number of full rows in 
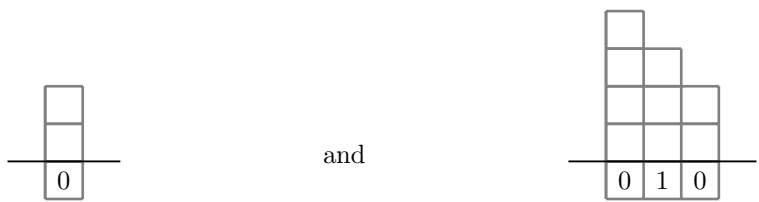
\begin{align*}
\begin{tikzpicture}[scale =.5] 
	\partitionfr{1,1,1}
	\draw [line width=0.25mm ] (-1,1) -- (2,1);
	\node at (.5,.5) {{$0$}};
\end{tikzpicture}
&&
\begin{tikzpicture}[scale =.5] 
	\node at (.5,.5) {{}};
	\node at (0,1.5) {{and}};
\end{tikzpicture}
&&
\begin{tikzpicture}[scale =.5] 
	\partitionfr{3,3,3,2,1}
	\draw [line width=0.25mm ] (-1,1) -- (4,1);
	\node at (.5,.5) {{$0$}};
	\node at (1.5,.5) {{$1$}};
	\node at (2.5,.5) {{$0$}};
\end{tikzpicture}
\end{align*}
are both $2$. Therefore, if we couldn't combine $\nu^1$ and $\nu^2$, then we cannot combine the pair
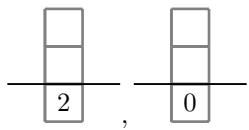
\begin{align*}
\begin{tikzpicture}[scale =.5] 
	\partitionfr{1,1,1}
	\draw [line width=0.25mm ] (-1,1) -- (2,1);
	\node at (.5,.5) {{$2$}};
\end{tikzpicture},
\begin{tikzpicture}[scale =.5] 
	\partitionfr{1,1,1}
	\draw [line width=0.25mm ] (-1,1) -- (2,1);
	\node at (.5,.5) {{$0$}};
\end{tikzpicture}
\end{align*}
\begin{proposition} We have the following:
\begin{enumerate}
\item Given a diagram with partition $\nu$ and image pair $\nu^1, \nu^2$, the number of full rows of $\nu$ is equal to the number of full rows of $\nu^1$. Thus, if $\nu^0 , \nu$ cannot be combined, then $\nu^0, \nu^1$ cannot be combined.
\item Suppose we have a sequence $\nu^1, \nu^2, \nu^3$ such that $\nu^1$ cannot combine with $\nu^2$, but $\nu^2$ can combine with $\nu^3$ to form $\nu$. Then $\nu^1$ cannot combine with $\nu$.
\end{enumerate}
\begin{proof}
\begin{enumerate}
\item This statement is obvious, since $\nu^1$ is the rightmost column of $\nu$. 
\item The number of full rows in $\nu$ equals the number of full rows of $\nu^2$, so that the inequality which determines whether we can combine $\nu^1$ and $\nu^2$ is still not satisfied between $\nu^1$ and $\nu$. 
\end{enumerate}
\end{proof}
\end{proposition}
Thus, we get the following sign-reversing involution:
Given $(\nu^1, \dots, \nu^{\ell(\mu)}) \in \Lt_{\mu}$ corresponding to the rearrangement $(\mu_{\alpha_1}, \dots, \mu_{\alpha_{\ell(\mu)}})$ proceed from $i=1$ to $i= \ell(\mu)$, until we find either the first $\nu^i$ which can combine with $\nu^{i+1}$, or the first $\nu^i$ for which $\mu_{\alpha_i} >1$. In the first case we combine $\nu^i$ and $\nu^{i+1}$. In the second case, we separate as described above. If no such $\nu^i$ is found, then leave the object fixed.

Our proposition above ensures that if a partition is separated, then reapplying the map will locate the pair which was created, since no preceding partitions were combined or separated. Likewise, if a pair was combined, then the map will locate the new partition which was formed. In either case the length of $\mu$ which underlies the sequence of partitions changes by exactly $1$. Therefore, the sign associated to $(\nu^1, \dots, \nu^{\ell(\mu)})$ will change parity. It remains to find the sequences which are fixed points:

Suppose $(\nu^1, \dots, \nu^{\ell(\mu)})$ corresponding to the rearrangement $(\mu_{\alpha_1}, \dots, \mu_{\alpha_{\ell(\mu)}})$ is a fixed point.
If $\mu_i>1$ for some $i$, then the object has a pair in the sign-reversing involution, so we must have $\mu_i =1 $ for all $i$. Therefore every fixed object is an element of $\Lt_{(1,\dots,1)}$ and has associated sign $(-1)^{(k+1)-(k+1)}= 1.$ This means that we are looking for columns of lengths
$(a_1, \dots, a_{k+1})$ and some underlying labels $(b_1,\dots, b_{k+1})$ respectively such that no two columns can be combined. In pictures, we have a sequence which looks like

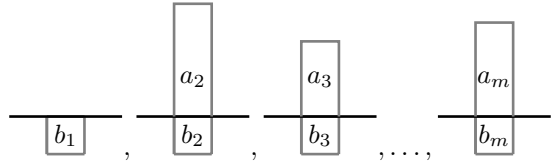
\begin{align*}
\begin{tikzpicture}[scale =.5] 
	\partitionfr{1}
	\draw [line width=1pt ] (-1,1) -- (2,1);
	\node at (.5,.5) {{$b_1$}};
\end{tikzpicture},
\begin{tikzpicture}[scale =.5] 
	\partitionfr{1}
	\draw [line width= 1pt ] (-1,1) -- (2,1);
	\draw [help lines,line width=1pt ] (0,1) -- (0,4) -- (1,4) -- (1,1);
	\node at (.5,.5) {{$b_2$}};
	\node at (.5,2) {{$a_2$}};
\end{tikzpicture},
\begin{tikzpicture}[scale =.5] 
	\partitionfr{1}
	\draw [line width=1pt ] (-1,1) -- (2,1);
	\draw [help lines,line width=1pt ] (0,1) -- (0,3) -- (1,3) -- (1,1);
	\node at (.5,.5) {{$b_3$}};
	\node at (.5,2) {{$a_3$}};
	\end{tikzpicture},
\dots,
\begin{tikzpicture}[scale =.5] 
	\partitionfr{1}
	\draw [line width=1pt ] (-1,1) -- (2,1);
	\draw [help lines, line width=1pt] (0,1) -- (0,3.5) -- (1,3.5) -- (1,1);
	\node at (.5,.5) {{$b_{m}$}};
	\node at (.5,2) {{$a_{m}$}};
\end{tikzpicture}
\end{align*}

We first see that $a_1 = 0$, since it is an element of $\Lt_{(1,\dots,1)}$. The condition that no two can be combined can be written as the inequality
$$
a_{i+1} < a_i +b_i.
$$
We therefore define a new set of objects $M^\lambda_k$ given by the set of all sequences 
$$
s = \left((a_1,b_1), (a_2,b_2), \dots, (a_{k+1},b_{k+1})\right)
$$
with the conditions that $a_1 = 0$; $(b_1,\dots, b_{k+1} )$ is a rearrangement of $\lambda_1,\dots, \lambda_{\ell(\lambda)}, 0^{k+1-\ell(\lambda)}$ (where there are $k+1-\ell(\lambda)$ zeroes added to $\lambda$); and $a_{i+1} < a_i + b_i$. For such a sequence $s$ define its weight by $$\rho(s) = a_1 + a_2+ \cdots + a_{k+1}.$$

This proves the main result:
\begin{theorem}
$$\langle \Delta_{e_k} e_n , f_{\lambda} \rangle \Big|_{q=1}= \sum_{s \in M^{\lambda}_k} t^{\rho(s)}. $$
\end{theorem}

\section{The bijection}
To prove the Delta Conjecture \cite{deltacon} at $q=1$, we must show 
\begin{theorem} 
$$
\sum_{s \in M^{\lambda}_k} t^{\rho(s)} = \sum_{\substack{D \in D_n \\ \lambda(D) = \lambda}} t^{\area(D)} H_{n-k}(D). 
$$
\begin{proof}
We first write the second sum as 
$$
 \sum_{\substack{D \in D_n \\ \lambda(D) = \lambda}} t^{\area({D})} H_{n-k}(D) = \sum_{\overline{D} \in D^{\lambda}_{n-k}} t^{\area(\overline{D})},
$$
where the objects in $D^{\lambda}_{n-k}$ are created in the following way:
\begin{enumerate}
\item Select $D \in D_n$ such that $\lambda(D) = \lambda$. Suppose its area sequence is $\alpha_1,\dots, \alpha_n$, and define the $0$ row with area $\alpha_0 = 0$ to be the origin. 
\item Select $n-k$ distinct rows $i_1,\dots, i_{n-k} \in \{0,2,\dots, n\}$ such that
either $i_j = 0$ or $\alpha_{i_j -1} = \alpha_{i_j}-1$. 
\item For all $j$, draw an $\times$ at each cell in row $i_j$ that contributes to the area of $D$. If $i_j = 0$, then mark the origin with an $\times$. After doing this for all $j$, call this new object $\overline{D}$.
\item Define the area of $\overline{D}$ by $\alpha_1+ \dots + \alpha_n - (\alpha_{i_1}+\cdots + \alpha_{i_{n-k}})$. That is, count all the cells with no $\times$ which contribute to the area of $D$.
\end{enumerate}
Here is an example: Suppose $D$ has area sequence $(0,1,2,2,1,2,1,2)$ so that the Dyck path is given by 
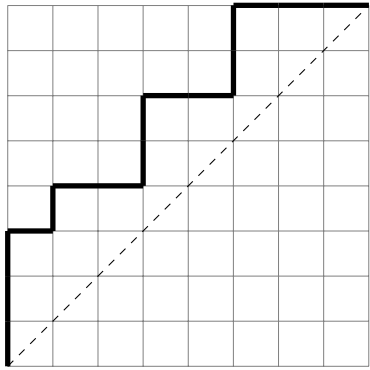
\begin{align*}
\begin{tikzpicture}[scale =.6] 
\dyckpath{0,0}{8}{1,1,1,0,1,0,0,1,1,0,0,1,1,0,0,0};
\end{tikzpicture}.
\end{align*}
We select rows $0,3,$ and $8$ to get the following labelled Dyck path in $D^{(3,2,2,1)}_{3}$:

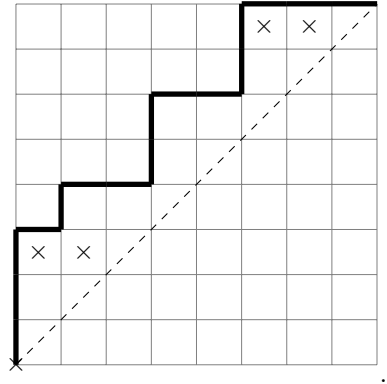
\begin{align*}
\begin{tikzpicture}[scale =.6] 
\dyckpath{0,0}{8}{1,1,1,0,1,0,0,1,1,0,0,1,1,0,0,0};
\node at (0,0) {{$\times$}};
\node at (.5,2.5) {{$\times$}};
\node at (1.5,2.5) {{$\times$}};
\node at (5.5,7.5) {{$\times$}};
\node at (6.5,7.5) {{$\times$}};
\end{tikzpicture}.
\end{align*}
Reading the rows bottom to top, the area of this object would be $0+1+0+2+1+2+1+0 = 7$. \\

The weight-preserving bijection $\phi:D^{\lambda}_{n-k} \rightarrow M^{\lambda}_k$ is given by the following steps:
\begin{enumerate}
\item To the left of every North step which begins a vertical segment say of length $b$ we write 
$(a,b)$, where $a$ is the area contribution of this row. 
\item To the left of every remaining North step whose row is not labelled, write $(c,0)$, where $c$ is the area contribution of this row.
\item For each label $(c,0)$, draw a North-East diagonal line from the beginning of its North step to the first start of an East step, and write $(c,0)$ above the East step which precedes it. (This East step exists since Dyck paths must always return to the main diagonal. So the only way the Dyck path can cross diagonal $c$ is with two consecutive East steps. This means every pair of consecutive North steps can be associated to a unique pair of consecutive East steps.)
\item Follow the Dyck path from bottom to top writing the pairs in the order which they appear.
\item If the origin is not labelled, include a $(0,0)$ at the end.
\end{enumerate}

It is easier to follow the bijection with pictures, so we illustrate these series of steps for a particular example. We select a Dyck path in $D^{(4,3,2,1)}_3$ and provide an element of $M^{(4,3,2,1)}_{7}$:\\
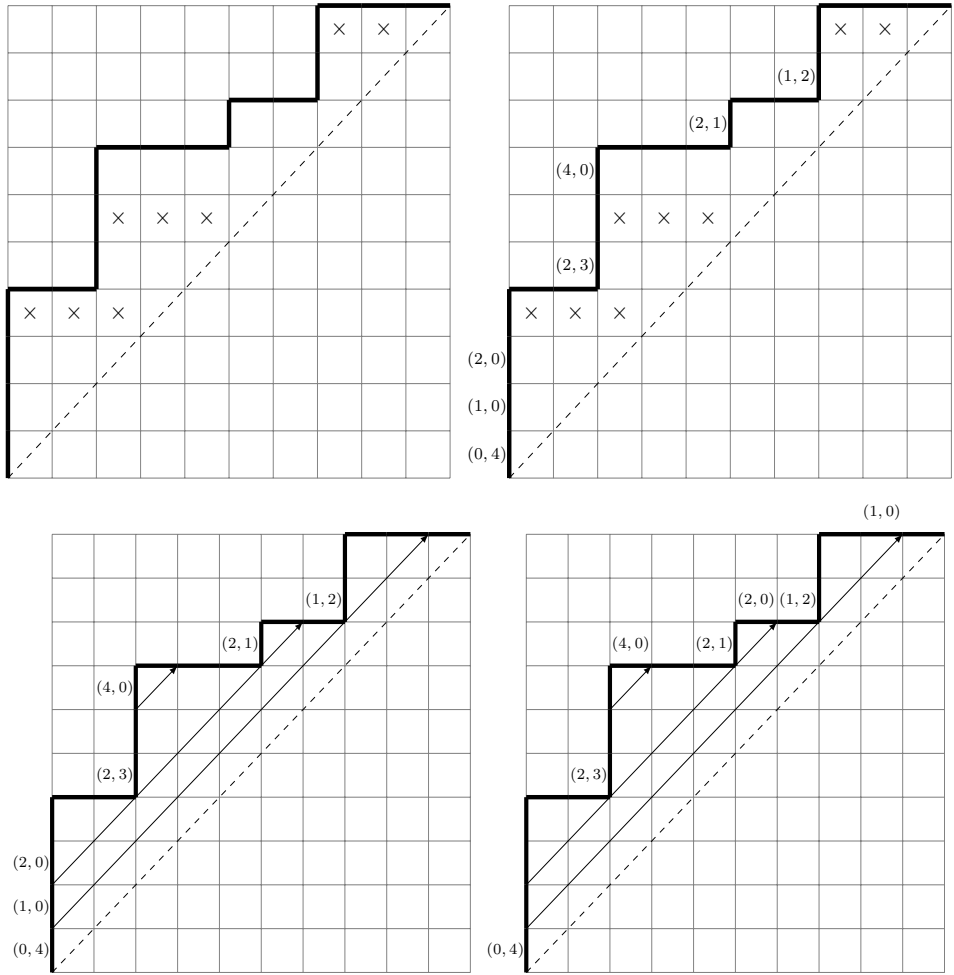
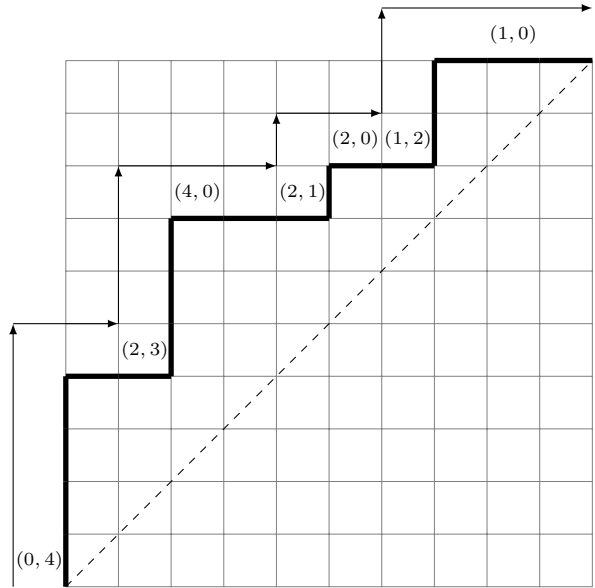
\begin{align*}
\resizebox{5in}{2.5in}{
\begin{tikzpicture}[scale =.7] 
\dyckpath{0,0}{10}{1,1,1,1,0,0,1,1,1,0,0,0,1,0,0,1,1,0,0,0};
\node at (.5,3.5) {{$\times$}};
\node at (1.5,3.5) {{$\times$}};
\node at (2.5,3.5) {{$\times$}};
\node at (2.5,5.5) {{$\times$}};
\node at (3.5,5.5) {{$\times$}};
\node at (4.5,5.5) {{$\times$}};
\node at (7.5,9.5) {{$\times$}};
\node at (8.5,9.5) {{$\times$}};
\end{tikzpicture}
\begin{tikzpicture}[scale =.7] 
\dyckpath{0,0}{10}{1,1,1,1,0,0,1,1,1,0,0,0,1,0,0,1,1,0,0,0};
\node at (.5,3.5) {{$\times$}};
\node at (1.5,3.5) {{$\times$}};
\node at (2.5,3.5) {{$\times$}};
\node at (2.5,5.5) {{$\times$}};
\node at (3.5,5.5) {{$\times$}};
\node at (4.5,5.5) {{$\times$}};
\node at (7.5,9.5) {{$\times$}};
\node at (8.5,9.5) {{$\times$}};
\node at (-.5,.5) {\scriptsize{$(0,4)$}};
\node at (-.5,1.5) {\scriptsize{$(1,0)$}};
\node at (-.5,2.5) {\scriptsize{$(2,0)$}};
\node at (1.5,4.5) {\scriptsize{$(2,3)$}};
\node at (1.5,6.5) {\scriptsize{$(4,0)$}};
\node at (4.5,7.5) {\scriptsize{$(2,1)$}};
\node at (6.5,8.5) {\scriptsize{$(1,2)$}};
\end{tikzpicture}}
\\
\resizebox{5in}{2.5in}{
\begin{tikzpicture}[scale =.7] 
\dyckpath{0,0}{10}{1,1,1,1,0,0,1,1,1,0,0,0,1,0,0,1,1,0,0,0};
\node at (-.5,.5) {\scriptsize{$(0,4)$}};
\node at (-.5,1.5) {\scriptsize{$(1,0)$}};
\node at (-.5,2.5) {\scriptsize{$(2,0)$}};
\node at (1.5,4.5) {\scriptsize{$(2,3)$}};
\node at (1.5,6.5) {\scriptsize{$(4,0)$}};
\node at (4.5,7.5) {\scriptsize{$(2,1)$}};
\node at (6.5,8.5) {\scriptsize{$(1,2)$}};
\draw[-{latex}] (0,1) -- (9,10);
\draw[-{latex}](0,2) -- (6,8);
\draw[-{latex}] (2,6) -- (3,7);
\end{tikzpicture}
\begin{tikzpicture}[scale =.7] 
\dyckpath{0,0}{10}{1,1,1,1,0,0,1,1,1,0,0,0,1,0,0,1,1,0,0,0};
\node at (-.5,.5) {\scriptsize{$(0,4)$}};
\node at (8.5,10.5) {\scriptsize{$(1,0)$}};
\node at (5.5,8.5) {\scriptsize{$(2,0)$}};
\node at (1.5,4.5) {\scriptsize{$(2,3)$}};
\node at (2.5,7.5) {\scriptsize{$(4,0)$}};
\node at (4.5,7.5) {\scriptsize{$(2,1)$}};
\node at (6.5,8.5) {\scriptsize{$(1,2)$}};
\draw[-{latex}][ ] (0,1) -- (9,10);
\draw[-{latex}](0,2) -- (6,8);
\draw[-{latex}] (2,6) -- (3,7);
\end{tikzpicture}
}
\end{align*}

\begin{align*}
\begin{tikzpicture}[scale =.7] 
\dyckpath{0,0}{10}{1,1,1,1,0,0,1,1,1,0,0,0,1,0,0,1,1,0,0,0};
\node at (-.5,.5) {\scriptsize{$(0,4)$}};
\node at (8.5,10.5) {\scriptsize{$(1,0)$}};
\node at (5.5,8.5) {\scriptsize{$(2,0)$}};
\node at (1.5,4.5) {\scriptsize{$(2,3)$}};
\node at (2.5,7.5) {\scriptsize{$(4,0)$}};
\node at (4.5,7.5) {\scriptsize{$(2,1)$}};
\node at (6.5,8.5) {\scriptsize{${(1,2)}$}};
\draw[-{latex}] (-1,0) -- (-1,5);
\draw[-{latex}](-1,5) -- (1,5);
\draw[-{latex}] (1,5) -- (1,8);
\draw[-{latex}] (1,8) -- (4,8);
\draw[-{latex}] (4,8) -- (4,9);
\draw[-{latex}] (4,9) -- (6,9);
\draw[-{latex}] (6,9) -- (6,11);
\draw[-{latex}] (6,11) -- (10,11);
\end{tikzpicture}
\end{align*}

Since the origin is not labeled, we include a $(0,0)$ at the end to get
$$
((0,4), (2,3),(4,0),(2,1),(2,0),(1,2),(1,0), (0,0)).
$$
For another example, we have
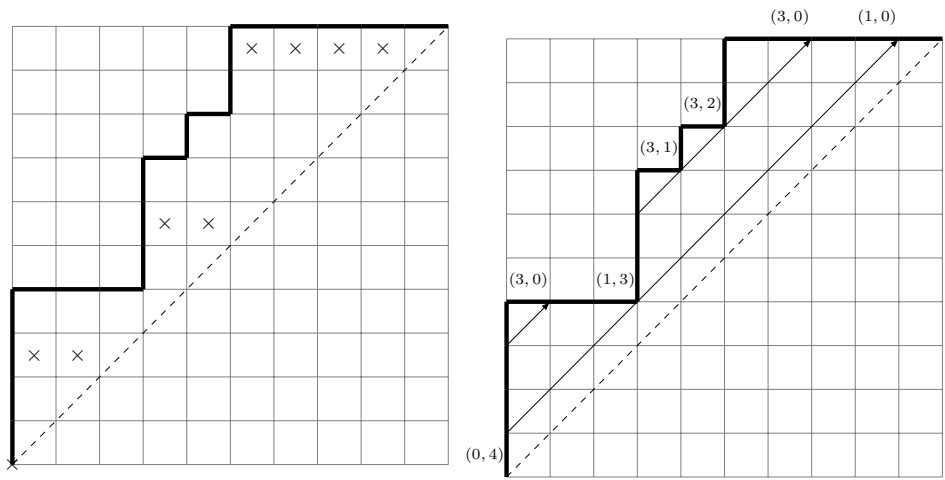
\begin{align*}
\resizebox{5in}{2.5in}{
\begin{tikzpicture}[scale =.7] 
\dyckpath{0,0}{10}{1,1,1,1,0,0,0,1,1,1,0,1,0,1,1,0,0,0,0,0};
\node at (0,0) {{$\bf\times$}};
\node at (.5,2.5) {{$\times$}};
\node at (1.5,2.5) {{$\times$}};
\node at (3.5,5.5) {{$\times$}};
\node at (4.5,5.5) {{$\times$}};
\node at (5.5,9.5) {{$\times$}};
\node at (6.5,9.5) {{$\times$}};
\node at (7.5,9.5) {{$\times$}};
\node at (8.5,9.5) {{$\times$}};
\end{tikzpicture}
\begin{tikzpicture}[scale =.7] 
\dyckpath{0,0}{10}{1,1,1,1,0,0,0,1,1,1,0,1,0,1,1,0,0,0,0,0};
\node at (-.5,.5) {\scriptsize{$(0,4)$}};
\node at (.5,4.5) {\scriptsize{$(3,0)$}};
\node at (2.5,4.5) {\scriptsize{$(1,3)$}};
\node at (3.5,7.5) {\scriptsize{$(3,1)$}};
\node at (4.5,8.5) {\scriptsize{$(3,2)$}};
\node at (6.5,10.5) {\scriptsize{$(3,0)$}};
\node at (8.5,10.5) {\scriptsize{$(1,0)$}};
\draw[-{latex}] (0,1) -- (9,10);
\draw[-{latex}](0,3) -- (1,4);
\draw[-{latex}] (3,6) -- (7,10);
\end{tikzpicture}}
\end{align*}

Since the origin is labeled, we do not include a $(0,0)$ at the end to get
$$
((0,4),(3,0),(1,3),(3,1),(3,2),(3,0),(1,0)).
$$

We must show that the image of $D^{\lambda}_{n-k}$ is in $M^{\lambda}_{k}$. It is easy to see that if $$\phi(\overline{D}) = ((a_1,b_1),\dots, (a_{k+1},b_{k+1})),$$ then by construction $(b_1,\dots,b_{k+1})$ is a rearrangement of
$\lambda, 0^{k+1-\ell(\lambda)}$. We also have $a_1 =0 $. It then suffices to show the inequalities are satisfied. Every two consecutive vertical line segments of lengths $b$ and $d$ (beginning at areas $a$ and $c$ respectively) will create a sequence in the image of the form 
$$
(a,b), (r_1,0),\dots, (r_{\ell},0),(c,d)
$$
This is described by the following diagram:
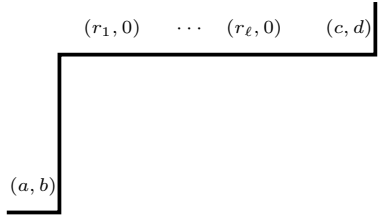
\begin{align*}
\begin{tikzpicture}[scale =.7] 
\node at (-.5,.5) {\scriptsize{$(a,b)$}};
\node at (1,3.5) {\scriptsize{$(r_1,0)$}};
\node at (2.5,3.5) {\scriptsize{$\cdots$}};
\node at (3.7,3.5) {\scriptsize{$(r_{\ell},0)$}};
\node at (5.5,3.5) {\scriptsize{$(c,d)$}};
\draw [line width=.5mm] (-1,0)--(0,0)-- (0,3) -- (6,3)-- (6,4);
\end{tikzpicture}
\end{align*}

We note that the top corner in the first vertical segment of the picture  lies on diagonal $a+b$. Likewise $r_1, \dots, r_\ell,c$ describe the right  diagonal of the east steps which they label. Therefore, since we are moving down diagonals as we proceed to the right, we have the inequalities
$$
a+b > r_1 > \dots > r_{\ell}>c
$$
This shows that $\phi(\overline{D})$ provides a sequence which satisfies the defining inequalities of $M^{\lambda}_k$and therefore $\phi$ is properly defined to give an element of $M^{\lambda}_k$. The area of the rows in $\overline{D}$ which contribute to the area are precisely the values of $a_1,\dots, a_{k+1}$. This means that $\area(\overline{D}) = \rho\left(\phi(\overline{D}) \right)$.
\\

To prove bijectivity, we need only provide the inverse. We do this first without proof, then show that this operation is indeed valid:
Given a sequence $((a_1,b_1),\dots, (a_{k+1},b_{k+1})) \in M^{\lambda}_k$, we do the following:
\begin{enumerate}
\item Proceeding from left to right, for each $(a,b)$ with $b \neq 0$ draw a line segment of length $b$ beginning on diagonal $a$, so that we form a Dyck path. 
\item For $(r_i,0)$ between $(a,b)$ and $(c,d)$, $b,d \neq 0$, draw $(r_i,0)$ over the east step ending in diagonal $r_i$. 
\item Move all $(r_i,0)$ South-West, down their diagonal, until we reach an end of a North step. Place $(r_i,0)$, to the left of the North step which lies above.
\item Label with $\times$'s the rows which have no pair to the left. 
\item If $(0,0)$ was not in the list, label the origin.
\end{enumerate}

We are simply doing the exact opposite from the definition of $\phi$. However, there are two implicit requirements which we have to check. The first step claims that the pairs $(a,b)$ with $b \neq 0$ when read from left to right describe a Dyck path. 
We again note that the sequence is composed of portions of the form
$$
(a,b), (r_1,0),\dots, (r_{\ell},0),(c,d).
$$
The fact that $c < r_\ell < \cdots < r_1 < a+b$ means that if we were to draw a vertical line segment of length $b$ such that the bottom North step contributes area $a$, then the diagonal of the top corner of the segment (this diagonal would be $a+b$) appears above diagonal $c$. So we have to place the segment of length $d$ strictly to the right of the segment of length $b$. Thus, the first step produces a Dyck path.
\\

The second place where some detail is needed is in step $2$, where we place $(r_i,0)$ between $(a,b)$ and $(c,d)$. Since $r_1> \cdots > r_{\ell}$, no two pairs will be placed above the same east step. Moreover, this portion of the Dyck path is horizontal starting from diagonal $a+b$ and ending on diagonal $c$, so there is indeed a unique East step which ends on diagonal $r_i$ for each $i$. This completes the bijection.
\end{proof}
\end{theorem}

\section{Further Comments}
Our computation reveals a rich collection of structures which can give the inner product of $\omega(\Delta_{e_k} e_n)$ in terms of any symmetric function F. 

\begin{theorem}

Let $F$ be a symmetric function of degree $n$ with expansion
\begin{align*}
F(x_1,\dots, x_{k+1}) &= c_1 x^{\alpha^1_1}_1 \cdots x^{\alpha^1_{k+1}}_{k+1} +  c_2 x^{\alpha^2_1}_1 \cdots x^{\alpha^2_{k+1}}_{k+1}+ \cdots +  c_r x^{\alpha^r_1}_1 \cdots x^{\alpha^r_{k+1}}_{k+1} \\
 &= \sum_i c_i  x^{\alpha^i_1}_1 \cdots x^{\alpha^i_{k+1}}_{k+1}.
\end{align*}
Then
$$
\langle \omega( \Delta_{e_k} e_n ) , F \rangle \Big|_{q=1} = 
\sum_{s \in L^F_k} {\weight(s)}
$$
where 
$L^F_k$ is the set of sequences
$$
s= ((a_1,b_1),(a_2, b_2), \dots, (a_{k+1}, b_{k+1})) 
$$
satisfying
\begin{enumerate}
\item $a_1 = 0$,
\item $(b_1,\dots, b_{k+1})  = (\alpha^j_1, \dots, \alpha^j_{k+1})$ for some $j$,
\item $a_{i+1} < a_i + b_i$, and
\item $\weight(s) = c_j t^{a_1+ \cdots + a_{k+1}}$.
\end{enumerate}
\end{theorem}

This theorem immediately gives that $\Delta_{e_k}e_n \Big|_{q=1}$ is t-positive in terms of the bases $\{s_{\lambda}\}_{\lambda}$, $\{m_{\lambda}\}_{\lambda}$, $\{e_{\lambda}\}_{\lambda}$, and surprisingly even $\{f_{\lambda}\}_{\lambda}$. We give here two particular examples which are important in studying the Delta Conjecture.

\subsection{The Hilbert series}

In recent work, Wilson \cite{tesler} was able to compute the Hilbert series $\langle \Delta_{e_k} e_n,p_1^n \rangle $ at $q=1$. His result can be obtained by our general method. We can write this inner-product as
$$
\langle \Delta_{e_k} e_n,p_1^n \rangle \Big|_{q=1}
= \sum_{\mu \vdash k+1} h_{\mu} \left[ \frac{1}{1-t} \right] f_\mu [1-t] p_{1}^n \left[\sum [\mu_i]_t \right].
$$

We can represent the symmetric function $p_1^n(x_1,\dots, x_{k+1})$ by 
$$
p_1^n(x_1,\dots, x_{k+1}) = \sum_{B_1,\dots, B_{k+1}} x_1^{|B_1|} \cdots x_{k+1}^{|B_{k+1}|}
$$
where $B_1,\dots, B_{k+1}$ are disjoint subsets of $\{1,\dots, n\}$ whose union is the whole set. In other words, we get an ordered set partition of length $k+1$, allowing empty sets. 

Let $P^n_k$ be the set of sequences
$$
P =((a_1, B_1),\dots, (a_{k+1}, B_{k+1}))
$$
of pairs where $a_1 = 0$; the $B_i$ are disjoint subsets of $\{1,\dots, n \}$ whose union is all of $\{1,\dots, n\}$; and 
$$
a_{i+1} < a_i + |B_i|.
$$
Letting $\rho(P) = a_1+\cdots + a_{k+1}$, we get
$$
\langle \Delta_{e_k} e_n,p_1^n \rangle \Big|_{q=1} = \sum_{P \in P^n_k} t^{\rho(P)}
$$
We can use the same bijection as in our previous section where the $b_i$ represent $|B_i|$, but now the sets describe cars along north segments. This means that the set of such sequences is equal to the set of parking functions in the $n \times n$ square with $n-k$ labeled rows, $t$-enumerated by the area statistic.  That is,

$$
\sum_{P \in P^n_k} t^{\rho(P)}
= \sum_{D\in D_n} H_{n-k}(D) \sum_{D(PF)=D} t^{\area(PF)},
$$
where $D(PF)$ is the supporting Dyck path of a parking function $PF$.

\subsection{A Schur function expansion}

The Delta Conjecture predicts that $\Delta_{e_k}e_n$ is $q,t$-Schur positive. We mean by this that the coefficient of $s_{\lambda}$ in $\Delta_{e_k}e_n$ is a polynomial in $q$ and $t$ with positive integer coefficients. 
Using the same methods we get the following expansion at $q=1$:

For a semistandard tableau $T$ let 
$$
c_i(T) = (\text{the number of times $i$ occurs in $T$}).
$$
Let $S^{\lambda}_{k}$ be the set of objects formed by 
\begin{enumerate}
\item selecting a semistandard filling $T$ of $\lambda$ with entries in $1,\dots, k+1$,
\item setting $a_1 = 0$ and choosing $a_2,\dots, a_{k+1}$ such that $a_{i+1} < a_i + c_i(T)$, and 
\item writing it all together as $S=( (a_1,\dots, a_{k+1}), T)$. 
\item Let $p(S) = a_1+ \cdots + a_{k+1}$.
\end{enumerate}

Then 
$$
\langle  \omega(\Delta_{e_k} e_n), s_{\lambda} \rangle \Big|_{q=1} = \sum_{S \in S^{\lambda}_k} t^{p(S)}.
$$
This gives the expansion
$$
\Delta_{e_k}e_n\Big|_{q=1} = \sum_{\lambda \vdash n} s_{\lambda}  \sum_{S \in S^{\lambda'}_k} t^{p(S)}.
$$
\section{Acknowledgements}
We must thank Adriano Garsia for his lessons on Macdonald polynomials and plethystic notation. Thank you for all the help and suggestions. 

\bibliographystyle{plain}
\bibliography{arXivDelta}

\begin{thebibliography}{10}

\bibitem{remarkable}
{A. M. Garsia and M. Haiman}.
\newblock {A remarkable $q,t$-Catalan sequence and $q$-Lagrange inversion}.
\newblock {\em J. Algebraic Combin.}, 5(3):191--244, 1996.

\bibitem{modifiedmac}
{A. M. Garsia and M. Haiman}.
\newblock {Some Natural Bigraded {$S_n$}-modules}.
\newblock {\em Electron. J. Combin.}, Volume 3(Issue 2 The Foata Festschrift),
  Jan 1996.

\bibitem{plethysm}
{A. M. Garsia, G. Xin, and M. Zabrocki}.
\newblock {Hall-Littlewood Operators in the Theory of Parking Functions and
  Diagonal Harmonics}.
\newblock {\em Int. Math. Res. Not.}, 6:1264--1299, Apr 2012.

\bibitem{ratshuffle}
{A. Mellit}.
\newblock {Toric braids and $(m,n)$-parking functions}.
\newblock {\em arXiv:1604.07456}, Apr 2016.

\bibitem{shuffle}
{E. Carlsson and A. Mellit}.
\newblock {A proof of the shuffle conjecture}.
\newblock {\em arXiv:1508.06239}, Aug 2015.

\bibitem{knots}
{E. Gorsky and A. Negu\c{t}}.
\newblock {Refined Knot Invariants and Hilber Schemes}.
\newblock {\em arXiv:1304.3328}, Apr 2013.

\bibitem{positivity}
{F. Bergeron, A. M. Garsia, M. Haiman, and G. Tesler}.
\newblock {Identities and Positivity Conjectures for some remarkable Operators
  in the Theory of Symmetric Functions}.
\newblock {\em Methods Appl. Anal.}, 6(3):363--420, Dec 1999.

\bibitem{sciencefiction}
{F. Bergeron and A. M. Garsia}.
\newblock {Science Fiction and Macdonald's Polynomials}.
\newblock {\em arXiv:9809128}, Sep 1998.

\bibitem{schroder}
J.~Haglund.
\newblock A proof of the $q,t$-{S}chr{\"o}der conjecture.
\newblock {\em Int. Math. Res. Not.}, 11:525--560, Aug 2004.

\bibitem{macdonald}
{I. G. Macdonald}.
\newblock {\em {Symmetric Functions and Hall Polynomials}}.
\newblock Oxford University Press, second edition, 1995.

\bibitem{deltacon}
{J. Haglund, J. B. Remmel, and A. T. Wilson}.
\newblock {The Delta Conjecture}.
\newblock {\em arXiv:1509.07058}, Sep 2015.

\bibitem{compos}
{J. Haglund, J. Morse, M. Zabrocki}.
\newblock A compositional shuffle conjecture specifying touch points of the
  dyck path.
\newblock {\em arXiv:1008.0828v1}, Aug 2010.

\bibitem{shufflecon}
{J. Haglund, M. Haiman, N. Loehr, J. B. Remmel, A. Ulyanov}.
\newblock {A Combinatorial Formula for the Character of the Diagonal
  Coinvariants}.
\newblock {\em arXiv:0310424}, Oct 2003.

\bibitem{n!}
{M. Haiman }.
\newblock {Hilbert schemes, polygraphs and the Macdonald positivity
  conjecture}.
\newblock {\em J. Amer. Math. Soc.}, 14:941--1006, May 2001.

\bibitem{bricktabloids}
{{\"O}. E{\~g}ecio{\~g}lu and J. B. Remmel}.
\newblock {Brick Tabloids and the connection matrices between bases of
  symmetric functions}.
\newblock {\em Discrete Appl. Math.}, 34:107--120, Nov 1991.

\bibitem{tewodros}
{T. Amdeberhan and E. Sergel Leven}.
\newblock {Multi-cores, posets, and lattice paths}.
\newblock {\em arXiv:1406.2250}, Jun 2015.

\bibitem{tesler}
A.~T. Wilson.
\newblock A weighted sum over generalized {Tesler} matrices.
\newblock {\em arXiv:1510.02684}, Oct 2015.

\end{thebibliography}

\end{document}